\input amstex
\documentstyle{amsppt}
\magnification=\magstep1
 \hsize 13cm \vsize 18.35cm \pageno=1
\loadbold \loadmsam
    \loadmsbm
    \UseAMSsymbols
\topmatter
\NoRunningHeads
\title An identity of the symmetry for the Frobenius-Euler polynomials
associated with the fermionic $p$-adic invariant $q$-integrals on
$\Bbb Z_p$
\endtitle
\author
  Taekyun Kim
\endauthor

\abstract The main purpose of this paper is to prove an identity of
symmetry for the Frobenius-Euler polynomials. It turns out that the
recurrence relation and multiplication theorem for the
Frobenius-Euler polynomials which discussed in [ K. Shiratani, S.
Yamamoto, On a $p$-adic interpolation function for the Euler numbers
and its derivatives, Memo. Fac. Sci. Kyushu University Ser.A,
39(1985), 113-125]. Finally we investigate several further
interesting properties of symmetry for the fermionic $p$-adic
invariant $q$-integral on $\Bbb Z_p$ associated with the
Frobenius-Euler polynomials and numbers.
\endabstract
\thanks  2000 AMS Subject Classification: 11B68, 11S80
\newline  This paper is supported by  Jangjeon Mathematical Society(JJMS-10R-2008)
\endthanks
\endtopmatter

\document

{\bf\centerline {\S 1. Introduction}}

 \vskip 20pt

The $n$-th Frobenius-Euler numbers  $H_n(q)$ and the $n$-th
Frobenius-Euler polynomials $H_n(q,x)$ attached to an algebraic
number $q (\neq 1)$ may be defined  by the exponential generating
functions
$$\sum_{n=1}^{\infty}H_n(q)\frac{t^n}{n!}=\frac{1-q}{e^t-q}, \text{
see [6,7], } \tag1$$
$$\sum_{n=0}^{\infty}H_n(q,x)\frac{t^n}{n!}=\frac{1-q}{e^t-q}e^{xt}.
$$
It is easy to show that $H_n(q,x)=\sum_{l=0}^n
\binom{n}{l}x^{n-l}H_l(q).$ Let $p$ be a fixed prime. Throughout
this paper $\Bbb Z_p ,$ $\Bbb Q_p ,$ $\Bbb C,$ and $\Bbb C_p$ will,
respectively, denote the ring of $p$-adic rational integers, the
field of $p$-adic rational numbers, the complex number field, and
the completion of algebraic closure of $\Bbb Q_p .$ When one talks
of $q$-extension, $q$ is variously considered as an indeterminate, a
complex $q \in \Bbb C$, or a $p$-adic number $q\in \Bbb C_p$, see
[9-22]. If $q\in \Bbb C$, then we assume $|q|<1.$ If $q\in \Bbb
C_p$, then we assume $|1-q|_p<1.$ For $x\in\Bbb Q_p ,$ we use the
notation $[x]_q=\frac{1-q^x}{1-q}, $ and $[x]_{-q}=\frac{1-(-q)^x
}{1+q},$ see [5-6]. The normalized valuation in $\Bbb C_p$ is
denoted by $|\cdot |_p$ with $|p|_p =\frac{1}{p} .$ We say that $f$
is a uniformly differentiable function at a point $a \in\Bbb Z_p $
and denote this property by $f\in UD(\Bbb Z_p )$, if the difference
quotients $F_f (x,y) = \dfrac{f(x) -f(y)}{x-y} $ have a limit
$l=f^\prime (a)$ as $(x,y) \to (a,a)$. For $f\in UD(\Bbb Z_p )$, let
us start with the expression
$$\eqalignno{ & \dfrac{1}{[p^N ]_q} \sum_{0\leq j < p^N} q^j f(j)
=\sum_{0\leq j < p^N} f(j)
\mu_q (j +p^N \Bbb Z_p ), }$$ representing a $q$-analogue of Riemann
sums for $f$, see [5, 6]. The integral of $f$ on $\Bbb Z_p$ will be
defined as limit ($n \to \infty$) of those sums, when it exists. The
$q$-deformed bosonic $p$-adic integral of the function $f\in UD(\Bbb
Z_p )$ is defined by
$$ I_q (f)= \int_{\Bbb Z_p }f(x) d\mu_q (x) = \lim_{N\to \infty}
\dfrac{1}{[dp^N ]_q} \sum_{0\leq x < dp^N} f(x) q^x, \text{ see [5]}
. $$
 Thus,  we note that
$$qI_{q}(f_1)= I_q(f)+(q-1)f(0)+\frac{q-1}{\log q}f^{\prime}(0),  $$
where $f_1(x)=f(x+1),$ $f^{\prime}(0)=\frac{df(0)}{dx}.$

The fermionic $p$-adic invariant $q$-integral on $\Bbb Z_p$ is
defined as
$$I_{-q}(f)=\int_{\Bbb Z_p}f(x)d\mu_{-q}(x)=\lim_{N\rightarrow
\infty}\frac{1}{[p^{N}]_{-q}}\sum_{x=0}^{p^N-1}f(x)(-q)^x, \text{
see [5]}.\tag2$$ In [8], H.J.H. Tuenter provided a generalization of
the Bernoulli number recurrence
$$B_m=\frac{1}{a(1-a^m)}\sum_{j=0}^{m-1}a^j\binom{m}{j}B_j\sum_{i=0}^{a-1}i^{m-j},
\text{ see [2, 3, 4]},$$ where $a, m \in \Bbb Z$ with $a>1$ $m\geq
1$, attributed to E.Y. Deeba and D.M. Rodriguez[2] and to I.
Gessel[3].  Define $S_m(k)=0^m+1^m+\cdots+k^m ,$ where $a, m \in
\Bbb Z$, with $a\geq 0$ and $m\geq 0 .$ H.J.H. Tuenter proved that
the quantity
$$\sum_{j=0}^m \binom{m}{j} a^{j-1}B_j b^{m-j}S_{m-j}(a-1), \text{
see [8], }$$ is symmetric in $a$ and $b$, provided $a, b, m \in \Bbb
Z,$ with $a>0, b>0 $ and $m\geq 0$. In this paper we prove  an
identity of symmetry for the Frobenius-Euler polynomials. It turns
out that the recurrence relation and multiplication theorem for the
Frobenius-Euler polynomials which discussed in [7]. Finally we
investigate the several further interesting properties of the
symmetry for the fermionic $p$-adic invariant $q$-integral on $\Bbb
Z_p$ associated with the Frobenius-Euler polynomials and numbers.

\vskip 20pt

{\bf\centerline {\S 2. An identity of symmetry for the
Frobenius-Euler polynomials}} \vskip 10pt

From (2) we can derive
$$qI_{-q}(f_1)+I_{-q}(f)=[2]_qf(0), \text{ where
$f_1(x)=f(x+1)$}.\tag3$$ By continuing this process, we see that
$$q^nI_{-q}(f_n)+(-1)^{n-1}I_{-q}(f)=[2]_q\sum_{l=0}^{n-1}(-1)^{n-1-l}q^lf(l),
\text{ where $f_n(x)=f(x+n)$.}$$ When $n$ is an odd positive
integer, we obtain
$$q^n I_{-q}(f_n)+I_{-q}(f)=[2]_q\sum_{l=0}^{n-1}(-1)^lf(l)q^l.
\tag4$$ If $n\in\Bbb N$ with $n\equiv 0 $ ($\mod 2$), then we have
$$q^n I_{-q}(f_n)-I_{-q}(f)=[2]_q\sum_{l=0}^{n-1}(-1)^{l-1}f(l)q^l
.\tag5$$ From (1) and (3) we derive
$$\int_{\Bbb Z_p}e^{xt}d\mu_{-q}(x)=\frac{1-
(-q)^{-1}}{e^t-(-q)^{-1}}=\sum_{n=0}^{\infty}H_n(-q^{-1})\frac{t^n}{n!}.
\tag6$$ Thus, we note that
$$\int_{\Bbb Z_p}x^n d\mu_{-q}(x)=H_n(-q^{-1}), \text{ and }
\int_{\Bbb Z_p}(y+x)^n d\mu_{-q}(x)=H_{n}(-q^{-1},x).$$

Let $n\in\Bbb N$ with $n\equiv 1$ ($\mod 2$). Then we obtain
$$[2]_q\sum_{l=0}^{n-1}(-1)^lq^l l^m =q^n H_{m}(-q^{-1},n)+
H_m(-q^{-1}).$$

For $ n\in \Bbb N$ with $n\equiv 0$ ($\mod 2$), we have
$$q^nH_m(-q^{-1},n)-H_m(-q^{-1})=[2]_q\sum_{l=0}^{n-1}(-1)^{l-1}q^l
l^m .$$

By substituting $f(x)=e^{xt}$ into (4), we can easily see that
$$\int_{\Bbb Z_p}q^ne^{(x+n)t} d\mu_{-q}(x)+\int_{\Bbb Z_p}e^{xt}
d\mu_{-q}(x)=[2]_{q}\frac{q^ne^{nt}+1}{qe^t+1}=[2]_q\sum_{l=0}^{n-1}(-1)^lq^le^{lt}.\tag7$$
Let $S_{k,q}(n)=\sum_{l=0}^{n}(-1)^ll^kq^k$. Then $S_{k,q}(n)$ is
called by the alternating sums of powers  of consecutive
$q$-integers. From the definition of the fermionic $p$-adic
invariant $q$-integral on $\Bbb Z_p$, we can derive

$$ \int_{\Bbb Z_p}q^n e^{(x+n)t} d\mu_{-q}(x)+\int_{\Bbb Z_p}e^{xt}
d\mu_{-q}(x)=\frac{[2]_q \int_{\Bbb Z_p} e^{xt}
d\mu_{-q}(x)}{\int_{\Bbb Z_p} e^{nxt}q^{(n-1)x}d\mu_{-q}(x)}.
\tag8$$

By (8), we easily see that
$$\int_{\Bbb Z_p}q^{(n-1)x}e^{nxt}
d\mu_{-q}(x)=\frac{1+q}{q^ne^{nt}+1}.$$

Let $w_1, w_2 (\in \Bbb N )$ be odd. By using double fermionic
$p$-adic invariant $q$-integral on $\Bbb Z_p$, we obtain
$$\frac{\int_{\Bbb Z_p }\int_{\Bbb Z_p}e^{(w_1x_1
+w_2x_2)t}d\mu_{-q}(x_1)d\mu_{-q}(x_2)}{\int_{\Bbb Z_p}
e^{w_1w_2xt}q^{(w_1w_2-1)x}d\mu_{-q}(x)}= \frac{[2]_{q} (q^{w_1
w_2}e^{w_1w_2t}+1)}{(qe^{w_1t}+1)(qe^{w_2t}+1)}.$$

Now we also consider the following fermionic  $p$-adic invariant
$q$-integral on $\Bbb Z_p$ associated with Frobenius-Euler
polynomials.

$$\aligned
I&=\frac{\int_{\Bbb Z_p }\int_{\Bbb Z_p}e^{(w_1x_1
+w_2x_2+w_1w_2x)t}d\mu_{-q}(x_1)d\mu_{-q}(x_2)}{\int_{\Bbb Z_p}
e^{w_1w_2xt}q^{(w_{1}w_{2}-1)x}d\mu_{-q}(x)}\\
&=\frac{[2]_q e^{w_1w_2xt}(q^{w_1w_2} e^{w_1w_2 t}+1
)}{(qe^{w_1t}+1)(qe^{w_2t}+1)}.
\endaligned\tag9$$

From (9) and (8), we can derive

$$\aligned
\frac{[2]_q\int_{\Bbb Z_p}e^{xt}d\mu_{-q}(x)}{\int_{\Bbb
Z_p}e^{w_1xt}q^{(w_1-1)x}d\mu_{-q}(x)}&=[2]_q \sum_{l=0}^{w_1-1}(-1)^lq^le^{lt}
=\sum_{k=0}^{\infty}\left([2]_q\sum_{l=0}^{w_1-1}(-1)^lq^ll^k\right)\frac{t^k}{k!}\\
&=\sum_{k=0}^{\infty}[2]_q S_{k,q}(w_1-1)\frac{t^k}{k!}.
\endaligned\tag10$$
By (9) and (10), we easily see that
$$\aligned
I&=\left(\frac{1}{[2]_q}\int_{\Bbb
Z_p}e^{w_1(x_1+w_2x)t}d\mu_{-q}(x)\right)\left(\frac{[2]_q\int_{\Bbb
Z_p}e^{w_{2} x_2t}d\mu_{-q}(x_2)}{\int_{\Bbb Z_p}e^{w_1 w_2
xt}q^{(w_1w_2-1)x}d\mu_{-q}(x)}\right)\\
&=\left(\frac{1}{[2]_q}\sum_{i=0}^{\infty}H_i(-q^{-1},
w_2x)\frac{w_1^i}{i!}t^i\right)\left([2]_q\sum_{l=0}^{\infty}S_{l,q^{w_2}}(w_1-1)\frac{w_2^l}{l!}t^l
\right)\\
&=\sum_{n=0}^{\infty}\left(\sum_{i=0}^n\binom{n}{i} H_i(-q^{-1},w_2
x)S_{n-i,q^{w_2}}(w_1-1)w_1^iw_{2}^{n-i}\right)\frac{t^n}{n!},
\endaligned\tag11$$
where $H_n(-q^{-1}, x)$ are the $n$-th Frobenius-Euler polynomials.

On the other hand,
$$\aligned
I&=\left(\frac{1}{[2]_q}\int_{\Bbb Z_p}e^{w_2(x_2+w_1
x)t}d\mu_{-q}(x_2)\right)\left(\frac{[2]_q\int_{\Bbb Z_p}e^{w_{1}
x_1t}d\mu_{-q}(x_1)}{\int_{\Bbb Z_p}e^{w_1 w_2
xt}q^{(w_1w_2-1)x}d\mu_{-q}(x)}\right) \\
&=\frac{1}{[2]_q}\left(\sum_{i=0}^{\infty}H_i(-q^{-1},w_1x)\frac{w_2^i
t^i}{i!}\right)\left([2]_q\sum_{l=0}^{\infty}S_{l,q^{w_1}}(w_2-1)\frac{w_1^l
t^l}{l!}\right)\\
&=\sum_{n=0}^{\infty}\left(\sum_{i=0}^n \binom{n}{i}H_i(-q^{-1}, w_1
x)S_{n-i, q^{w_1}}(w_2-1)w_2^iw_1^{n-i}\right)\frac{t^n}{n!}.
\endaligned\tag12$$
By comparing the coefficients on the both sides of (11) and (12), we
obtain the following theorem.

\proclaim{ Theorem 1} Let $w_1, w_2 (\in \Bbb N)$ be odd and let
$n(\geq 0)$ with $n \equiv 1 (\mod 2)$. Then we have
$$\aligned
&\sum_{i=0}^n \binom{n}{i}H_i(-q^{-1},w_2 x)
S_{n-i,q^{w_2}}(w_1-1)w_1^iw_2^{n-i}\\
&=\sum_{i=0}^n
\binom{n}{i}H_i(-q^{-1},w_1x)S_{n-i,q^{w_1}}(w_2-1)w_2^iw_1^{n-i},
\endaligned \tag13$$
where $H_n(q,x)$ are the $n$-th Frobenius-Euler polynomials.
\endproclaim
Setting $x=0$ in (13), we obtain the following corollary.

\proclaim{ Corollary 2} Let $w_1, w_2 (\in \Bbb N)$ be odd and let
$n\in\Bbb Z_{+}$ be an odd. Then we have
$$\aligned
&\sum_{i=0}^n \binom{n}{i}H_i(-q^{-1})S_{n-i,q^{w_2}}(w_1-1)w_1^i
w_2^{n-i}\\
&=\sum_{i=0}^n
\binom{n}{i}H_i(-q^{-1})S_{n-i,q^{w_1}}(w_2-1)w_2^iw_1^{n-i},
\endaligned$$
where $H_{i}(-q^{-1})$ are the $n$-th Frobenius-Euler numbers.
\endproclaim

If we take $w_2=1$ in (13), then we have
$$H_n(-q^{-1},w_1x)=\sum_{i=0}^n\binom{n}{i}H_i(-q^{-1},x)S_{n-i,q}(w_1-1)w_1^i.\tag14$$
Setting $x=0$ in (14), we obtain the following corollary.

\proclaim{ Corollary 3} Let $w_1 (>1)$ be an odd integer and let
$n\in\Bbb Z_{+}$ with $n\equiv 1 (\mod 2)$. Then we have
$$H_n(-q^{-1})=\frac{1}{1-w_1^n}\sum_{i=0}^{n-1}\binom{n}{i}H_i(-q^{-1})S_{n-i,q}(w_1-1)w_1^i.$$
\endproclaim
From (7) and (8), we derive
$$\aligned
I&=\left(\frac{e^{w_1w_2xt}}{[2]_q}\int_{\Bbb
Z_p}e^{w_1x_1t}d\mu_{-q}(x_1)\right)\left(\frac{[2]_q\int_{\Bbb
Z_p}e^{w_2x_2t}d\mu_{-q}(x_2)}{\int_{\Bbb
Z_p}e^{w_1w_2xt}q^{(w_1w_2-1)x}d\mu_{-q}(x)}\right)\\
&=\left(\frac{e^{w_1w_2xt}}{[2]_q}\int_{\Bbb
Z_p}e^{w_1x_1t}d\mu_{-q}(x_{1})\right)\left([2]_q\sum_{l=0}^{w_1-1}(-1)^lq^{w_2
l}e^{w_2 lt}\right)\\
&=\sum_{l=0}^{w_1-1}(-1)^l q^{w_2l}\int_{\Bbb
Z_p}e^{(x_1+w_2x+(\frac{w_2}{w_1})l)tw_1}d\mu_{-q}(x_1)\\
&=\sum_{n=0}^{\infty}\left(w_1^n\sum_{l=0}^{w_1-1}(-1)^lq^{w_2l}
H_n(-q^{-1},w_2x+\frac{w_2}{w_1}l )\right)\frac{t^n}{n!}.
\endaligned\tag15$$
On the other hand,
$$\aligned
I&=\left(\frac{e^{w_1w_2xt}}{[2]_q}\int_{\Bbb
Z_p}e^{w_2x_2t}d\mu_{-q}(x_2)\right)\left(\frac{[2]_q\int_{\Bbb
Z_p}e^{w_1x_1t}d\mu_{-q}(x_1)}{\int_{\Bbb
Z_p}e^{w_1w_2xt}q^{(w_1w_2-1)x}d\mu_{-q}(x)}\right)\\
&=\left(\frac{1}{[2]_q}\int_{\Bbb
Z_p}e^{w_2x_2t}d\mu_{-q}(x_{2})\right)\left([2]_q\sum_{l=0}^{w_2-1}(-1)^lq^{w_1
l}e^{(w_1 l +w_1w_2 x)t}\right)\\
&=\sum_{l=0}^{w_2-1}(-1)^l q^{w_1l}\int_{\Bbb
Z_p}e^{(x_2+w_1x+\frac{w_1}{w_2}l)tw_2}d\mu_{-q}(x_2)\\
&=\sum_{n=0}^{\infty}\left(w_2^n\sum_{l=0}^{w_2-1}(-1)^lq^{w_1l}
H_n(-q^{-1},w_1x+\frac{w_1}{w_2}l )\right)\frac{t^n}{n!}.
\endaligned\tag16$$
By comparing the coefficients on the both sides of 915) and (160, we
obtain the following theorem.

\proclaim{ Theorem 4} Let $w_1, w_2 (\in \Bbb N)$ be odd and let
$n\in \Bbb Z_{+}$ with $n \equiv 1 (\mod 2)$. Then we have
$$w_1^n\sum_{l=0}^{w_1-1}(-1)^lq^{w_2l}H_n(-q^{-1},
w_2x+\frac{w_2}{w_1}l)
=w_2^n\sum_{l=0}^{w_2-1}(-1)^lq^{w_1l}H_n(-q^{-1},
w_1x+\frac{w_1}{w_2}l). $$
\endproclaim

Setting $w_2=1$ in Theorem 4, we get the multiplication theorem for
the Frobenius-Euler polynomials as follows:
$$H_n(-q^{-1},w_1x)=w_1^n\sum_{l=0}^{w_1-1}(-1)^lq^l
H_n(-q^{-1},x+\frac{l}{w_1}).$$

\vskip 20pt

\quad Taekyun Kim

\quad Division of General-Education, Kwangwoon University, Seoul
139-701, S. Korea

\quad e-mail:\text{ tkim$\@$kw.ac.kr; tkim64$\@$hanmail.net}

\enddocument